\input amstex\documentstyle{amsppt}  
\pagewidth{12.5cm}\pageheight{19cm}\magnification\magstep1
\topmatter
\title On the definition of almost characters\endtitle
\author G. Lusztig\endauthor
\address{Department of Mathematics, M.I.T., Cambridge, MA 02139}\endaddress
\thanks{Supported in part by National Science Foundation grant 1303060.}\endthanks  \endtopmatter   
\document

\define\mpb{\medpagebreak}

\define\si{\sim}

\define\sqc{\sqcup}

\define\qua{\quad}

\define\bA{\bar A}

\define\bZ{\bar Z}
\define\lb{\linebreak}

\define\part{\partial}
\define\emp{\emptyset}

\define\ra{\rangle}
\define\n{\notin}

\define\m{\mapsto}

\define\la{\langle}

\define\sub{\subset}    

\define\T{\times}
\define\ti{\tilde}
\define\nl{\newline}
\redefine\i{^{-1}}

\define\un{\underline}

\define\ot{\otimes}

\define\tr{\text{\rm tr}}

\define\a{\alpha}
\redefine\b{\beta}
\redefine\c{\chi}

\redefine\d{\delta}
\define\e{\epsilon}
\define\et{\eta}
\define\io{\iota}

\define\ph{\phi}

\define\r{\rho}
\define\s{\sigma}
\redefine\t{\tau}

\redefine\l{\lambda}
\define\z{\zeta}
\define\x{\xi}

\redefine\G{\Gamma}
\redefine\D{\Delta}

\define\Ps{\Psi}

\define\bb{\bold b}

\define\kk{\bold k}

\define\CC{\bold C}

\define\FF{\bold F}

\define\ZZ{\bold Z}

\define\cb{\Cal B}
\define\cc{\Cal C}

\define\ce{\Cal E}
\define\cf{\Cal F}
\define\cg{\Cal G}
\define\ch{\Cal H}

\define\co{\Cal O}

\define\car{\Cal R}

\define\fs{\frak s}

\define\tch{\ti\ch}

\define\ty{\ti y}

\define\tC{\ti C}

\define\DL{DL}
\define\EEIGHT{L1}
\define\ORA{L2}
\define\ICM{L3}
\define\FAM{L4}

\head Introduction\endhead
\subhead 0.1\endsubhead
Let $\cg(\FF_q)$ be the group of rational points of a connected reductive group 
$\cg$ defined over a finite field $\FF_q$. The vector space of $\CC$-valued class 
functions on $\cg(\FF_q)$ has a basis $\b$ given by the characters of the
irreducible representations of $\cg(\FF_q)$. It has a second basis $\b'$
(whose elements are defined only up to multiplication by roots of $1$) given by the characteristic functions of the Frobenius-invariant character sheaves on $\cg$. 
To compare these two bases it is convenient to introduce a third basis $\b''$
(whose elements are again defined only up to multiplication by roots of $1$). Now 
the functions in $\b''$ are given by explicit linear combinations of elements of 
$\b$;  they are called almost characters. In this paper we give a definition of 
$\b''$. (A definition of almost characters assuming that $\cg$ has connected 
centre appeared in \cite{\ORA}.) The interest of the basis
$\b''$ is that it conjecturally coincides with $\b'$. 

This paper was written in response to a question that Meinolf Geck asked me; I 
wish to thank him for asking the question and for comments on an earlier version.

\subhead 0.2. Notation\endsubhead
$\kk$ denotes a fixed algebraically closed field of characteristic $p\ge0$.
All algebraic groups are assume to be affine and over $\kk$. 
If $H$ is an algebraic group we denote by $H^0$ the identity component of $H$. 
For $g\in H$ we denote by $g_s$ (resp. $g_u$) the semisimple (resp. unipotent) part of $g$, so that $g=g_sg_u=g_ug_s$. 

If $x$ is an element of a group $\G$ we denote by $Z_\G(x)$ the centralizer of $x$ 
in $\G$.

\head 1. Nonabelian Fourier transform\endhead
\subhead 1.1\endsubhead
In this paper a representation of a group $\G$ is always assumed to be in
a finite dimensional $\CC$-vector space. 

We now assume that $\G$ is a finite group. Let $M(\G)$ be the set 
consisting of all pairs $(x,\s)$ where $x\in\G$ and $\s$ is an irreducible
representation of $Z_\G(x)$ (up to isomorphism); these pairs are taken 
modulo the equivalence relation $(x,\s)\m(vxv\i,{}^v\s)$ for any
$v\in\G$ (here ${}^v\s$ is the irreducible representation of $Z_\G(vxv\i)$
defined by composing $\s$ with $z\m v\i zv$). Following \cite{\EEIGHT} we
define a pairing $\{,\}:M(\G)\T M(\G)@>>>\CC$ by
$$\align&\{(x,\s),(y,\t)\}\\&=\sum_{w\in\G;xwyw\i=wyw\i x}|Z_\G(x)|\i||Z_\G(y)|\i
\tr(w\i x\i w,\t)\tr(wyw\i,\s).\tag a\endalign$$

\subhead 1.2\endsubhead
Let $c$ be an integer $\ge1$. Let $\mu_c=\{\a\in\CC;\a^c=1\}$.

A $\ZZ/c$-grading for a finite group $\G$ is
by definition a surjective homomorphism $\z:\G@>>>\ZZ/c$ 
such that for some $x\in\z\i(1)$ we have $x^c=1$. For such $\G,\z$ and for $i=0,1$
let $M^i_\z(\G)$ be the subset of $M(\G)$ consisting of all $(x,\s)\in M(\G)$ 
such that $x\in\z\i(i)$ and

(a) $Z_\G(x)\cap\z\i(1)\ne\emp$, $\s|_{Z_\G(x)\cap\z\i(0)}$ is irreducible.
\nl
(When $i=1$ the conditions (a) are automatically satisfied.)

For $\a\in\mu_c$ let $\e_\a:\G@>>>\CC^*$ be the homomorphism which takes 
any $x\in\z\i(1)$ to $\a$; the restriction of $\e_\a$ to various subgroups
of $\G$ is denoted again by $\e_\a$.

Let $V(\G)$ be the vector space of formal linear combinations 
$\sum_{(x,\s)\in M(\G)}b_{x,\s}(x,\s)$ with $b_{x,\s}\in\CC$. For $i=0,1$
we define a subspace $V^i_\z(\G)$ of $V(\G)$ by the conditions

$b_{x,\s}=0$ if $(x,\s)\n M^i_\z(\G)$,

$b_{x,\s\ot\e_\a}=\a^{i-1}b_{x,\s}$ for any $\a\in\mu_c$.
\nl
Let $A:V(\G)@>>>V(\G)$ be the linear map given by
$$A(x,\s)=\sum_{(y,\t)\in M(\G)}\{(x,\s),(y,\t)\}(y,\t)$$
for any $(x,\s)\in M(\G)$. For $i=0,1$ we show that 

(b) $A(V^i_\z(\G))\sub V^{1-i}_\z(\G)$.
\nl
Since $V^i_\z(\G)$ is generated by the vectors
$\sum_{\a\in\mu_c}\a^{i-1}(x,\s\ot\e_\a)$ for various 
$(x,\s)\in M^i_\z(\G)$, it is enough to show that for $(x,\s)\in M^i_\z(\G)$ 
we have 

$A(\sum_{\a\in\mu_c}\a^{i-1}(x,\s\ot\e_\a))\in V^{1-i}_\z(\G)$ or equivalently

(c) $\sum_{(z,\r)\in M(\G)}b'_{z,\r}(z,\r)\in V^{1-i}_\z(\G)$
\nl
where 
$$\align&b'_{z,\r}=\sum_{w\in\G;xwzw\i=wzw\i x}|Z_\G(x)|\i|Z_\G(z)|\i\\&
\T\tr(w\i x\i w,\r)\sum_{\a\in\mu_c}\a^{i-1}\tr(wzw\i,\s\ot\e_\a).\endalign$$
The last sum is equal to $c\d_{1-i,\z(z)}\tr(wzw\i,\s)$
hence, if $b'_{z,\r}\ne0$, then $z\in\z\i(1-i)$.
If in addition we have $i=1$ then $Z_\G(z)$ contains
$w\i xw\in\z\i(1)$ (for some $w$) which belongs to $Z_\G(z)$.
Thus $Z_\G(z)\cap\z\i(1)\ne0$; moreover in this case we have
$\tr(w\i x\i w,\r)\ne0$ for some $w$ and $w\i x\i w\in\z\i(-1)$; 
this forces $\r|_{Z_\G(z)\cap\z\i(0)}$ to be irreducible.
 
For $i=0,1$ we have 
$\tr(w\i x\i w,\r\ot\e_{\a'})=\a'{}^{-i}\tr(w\i xw,\r)$ for any 
$\a'\in\mu_c$, so that $b'_{z,\r\ot\e_{\a'}}=\a'{}^{-i}b'_{z,\r}$ for
any $\a',z,\r$. We see that (c) holds and (b) is proved.
From \cite{\EEIGHT} it is known that $A^2=1$. Combining this with (b) we 
see that

(d) {\it $A:V^1_\z(\G)@>>>V^0_\z(\G)$, $A:V^0_\z(\G)@>>>V^1_\z(\G)$ are 
inverse isomorphisms.}
\nl
(This was stated without proof in \cite{\ORA, (4.16.1)}).

\subhead 1.3\endsubhead
A $\ZZ$-grading for a group $\D$ is by definition a surjective homomorphism
$\et:\D@>>>\ZZ$ such that $\et\i(0)$ is finite. We now fix such $\D,\et$.
Let $\D^*$ be the set of central elements $\x\in\D$ such that $\x=y^c$ for
some $y\in\et\i(1)$ and some $c\ge1$. We have $\D^*\ne\emp$; indeed, if
$y\in\et\i(1)$, then $y^c$ is in the centre of $\D$ for some $c\ge1$. Moreover,

(a) {\it we can find $d\ge1$ such that $y\m y^d$ is a constant function
on $\et\i(1)$; for this $d$ we have $y^d\in\D^*$ for all $y\in\et\i(1)$.}
\nl
If $\x\in\D^*,c'\ge1$ then $\x^{c'}\in\D^*$. Moreover,

(b) {\it if $\x\in\D^*,\x'\in\D^*$ then $\x^a=\x'{}^{a'}$ for some $a\ge1,a'\ge1$.}
\nl
Indeed, we have $\x=y^c$, $\x'=y'{}^{c'}$ for some $c\ge,c'\ge1$ and some 
$y,y'$ in $\et\i(1)$. Hence if $d$ is as in (a), we have 
$\x^{c'd}=y^{cc'd}=y'{}^{cc'd}=\x'{}^{cd}$, so that (b) holds with $a=c'd,a'=cd$.

\mpb

For $i=0,1$ let $M^i_\et(\D)$ be the set of all pairs $(x,\s)$ where
$x\in\et\i(i)$ and $\s$ (an irreducible representation up to isomorphism
of $Z_\D(x)\cap\et\i(0)$) are such that

(c) {\it $Z_\D(x)\cap\et\i(1)\ne\emp$, $\s$ can be extended to a representation
of $Z_\D(x)$;}
\nl
these pairs are taken modulo the equivalence relation $(x,\s)\si(vxv\i,{}^v\s)$
for any $v\in\D$ (here ${}^v\s$ is the irreducible representation of
$Z_\D(vxv\i)\cap\d\i(0)$ defined by composing $\s$ with $z\m v\i zv$).

(When $i=1$ the conditions (c) are automatically satisfied.)

For any $\x\in\D^*$ we set $\D_\x=\D/\la\x\ra$ where $\la\x\ra$ is the 
cyclic subgroup of the centre of $\D$ generated by $\x$; this is a finite group.
For any $x\in\D$ we denote by $x_\x$ the image of $x$ in $\D_\x$. Now
$\et$ induces a surjective homomorphism $\z:\D_\x@>>>\ZZ/\et(\x)$ which is
a $\ZZ/\et(\x)$-grading for $\D_\x$. (We can find $y\in\et\i(1)$ such that
$y^{\et(\x)}=\x$; then $y_\x^{\et(\x)}=1$ and $\z(y_\x)=1$.)

For $i=0,1$ let $U^i_\et(\D)$ be the vector space of formal linear 
combinations $\sum_{(y,\t)\in M^i_\et(\D)}b_{y,\t}(y,\t)$ with 
$b_{x,\s}\in\CC$. To any $(x,\s)\in M^0_\et(\D)$ we associate an element 
$(x,\s)\hat{}\in U^1_\et(\D)$ (defined up to multiplication by a root of 
$1$) as follows.

We choose $\x\in\D^*$ and an extension of $\s$ to a
representation $\ti\s$ of $Z_\G(x)$ on which $\x$ acts as $1$; we set
$$\align&(x,\s)\hat{}_{\x,\ti\s}\\&=\sum_{(y,\t)\in M^1_\et(\D)}
\sum_{w\in\D/\la\x\ra;xwyw\i=wyw\i x}\et(\x)
|Z_\D(x)/\la\x\ra|\i|Z_\D(y)/\la\x\ra|\i\\&
\T\tr(w\i x\i w,\ti\t)\tr(wyw\i,\ti\s)(y,\t)\in U^1_\et(\D)\tag d\endalign$$
where $\ti\t$ is the representation of $Z_\D(y)$ extending $\t$ on which
$y$ acts as $1$. 

Assume now that $\x,\ti\s$ are replaced by another pair $\x',\ti\s'$ with 
the same properties as $\x,\ti\s$. We want to show that

(e) $(x,\s)\hat{}_{\x',\ti\s'}$ equals $(x,\s)\hat{}_{\x,\ti\s}$
times a root of $1$.
\nl
Using (b) we can assume that $\x'=\x^m$ for some $m\ge1$. From the 
definitions we see that 

$(x,\s)\hat{}_{\x^m,\ti\s}=(x,\s)\hat{}_{\x,\ti\s}$
\nl
and 

$(x,\s)\hat{}_{\x^m,\ti\s}$ is equal to $(x,\s)\hat{}_{\x^m,\ti\s'}$ times a root of $1$
\nl
so that (e) follows. We shall write $(x,\s)\hat{}$ instead of
$(x,\s)\hat{}_{\x,\ti\s}\in U^1_\et(\D)$; this is 
well defined only up to multiplication by a root of $1$.

\subhead 1.4\endsubhead
For any $(x,\s)\in M^0_\et(\D)$ we choose $\x\in\D^*$ and an extension of 
$\s$ to a representation $\ti\s$ of $Z_\G(x)$ on which $\x$ acts as $1$; we
can assume that $\x$ is independent of $(x,\s)$. We show:

(a) {\it The elements $(x,\s)\hat{}_{\x,\ti\s}$ for various
$(x,\s)\in M^0_\et(\D)$ form a $\CC$-basis of $U^1_\et(\D)$.}
\nl
We set $c=\et(\x)$. Let $\z:\D_\x@>>>\ZZ/c$ be as in 1.3. We define a 
linear map $f:U^0_\et(\D)@>>>V^0_\z(\D_\x)$ by 

$(x,\s)\m\sum_{\a\in\mu_c}\a\i(x_\x,\ti\s\ot\e_\a)$
\nl
(here we view $\ti\s$ as an irreducible representation of 
$Z_{\D_\x}(x_\x)$); this is easily seen to be an isomorphism.
We define a linear map $f':V^1_\z(\D_\x)@>>>U^1_\et(\D)$ by
$$\sum_{(y,\t)\in M^1_\z(\D_x)}b_{y,\t}(y,\t)\m
\sum_{(y,\t)\in M^1_\z(\D_\x); y=1\text{ on }\t}b_{y,\t}
(\ty,\t|_{Z_{\D_\x}(y)\cap\z\i(0)})$$
where for $(y,\t)\in M^1_\z(\D_\x)$ we define $\ty\in\et\i(1)$ by 
$\ty_\x=y$. (Note that $Z_{\D_\x}(y)\cap\z\i(0)=Z_\D(\ty)\cap\et\i(0)$.)
This is easily seen to be an isomorphism. The composition
$$V^0_\D@>f>>V^0_{\D_\x}@>A>>V^1_{\D_\x}@>f'>>V^1_\D$$
is an isomorphism (see 1.2(d)). It takes $(x,\s)\in M^0_\et(\D)$ to 
$(x,\s)\hat{}_{\x,\ti\s}$. This proves (a).

\head 2. Special elements\endhead
\subhead 2.1\endsubhead
In this section we assume that $p=0$. We fix an integer $q\ne0$.

Let $\ch$ be a connected reductive group.
For any element $h\in\ch$ let $\cb_u$ be the variety of Borel subgroups of 
$\ch$ that contain $h$.

Let $\ch_u^{sp}$ be the set of special unipotent elements of $\ch$, in
the sense of \cite{\ORA, (13.1.1)}; for $u\in\ch_u^{sp}$ let
$\bA(u)$ be the canonical quotient group of $Z_{\ch}(u)/Z^0_{\ch}(u)$
defined in  \cite{\ORA, 13.1} (paragraph following (13.1.2)) and let
$I_{\ch}(u)$ be the kernel of the obvious homomorphism
$Z_{\ch}(u)/Z^0_{\ch}(u)@>>>\bA(u)$. (When we refer to \cite{\ORA, 13.1}, we
assume that the ground field in {\it loc. cit.} is replaced by our 
$\kk$.) For further information on $A(u)$ see \cite{\FAM}. We show:

(a) {\it Assume that $u\in\ch^{sp}_u$. Then $u^q\in\ch^{sp}_u$ and we have 
canonically 

$Z_{\ch}(u)/Z_{\ch}^0(u)=Z_{\ch}(u^q)/Z_{\ch}^0(u^q)$, $I_{\ch}(u)=I_{\ch}(u^q)$.}
\nl
We have clearly $Z_{\ch}(u)\sub Z_{\ch}(u^q)$. Since $q\ne0$ and $p=0$, $u$ 
is conjugate to $u^q$ in $H$; it follows that
$Z_{\ch}(u)=Z_{\ch}(u^q)$. This implies the first equality in (a). We have clearly
$\cb_u\sub\cb_{u^q}$. Since $\cb_u,\cb_{u^q}$ are isomorphic, we must have
$\cb_u=\cb_{u^q}$. Using this and the definition we see that the second 
equality in (a) holds.

\subhead 2.2\endsubhead
In the remainder of this paper $G$ denotes a connected reductive group over 
$\kk$. Moreover we fix an automorphism $j:G@>\si>>G$ of finite order.

\subhead 2.3\endsubhead
We say that $g\in G$ is special if 

(i) $g_s^{q^r}=g_s$ for some $r\ge1$ and 

(ii) $g_u$ is a special unipotent element of the connected reductive group $Z^0_G(g_s)$.
\nl
Let $G^{sp}$ be the set of special elements of $G$. 
For $g\in G^{sp}$, setting $\tch=Z_G(g_s)$, $\ch=\tch^0$, we have obvious imbeddings 
$$I_{\ch}(g_u)\sub Z_{\ch}(g_u)/Z^0_{\ch}(g_u)\sub Z_{\tch}(g_u)/Z^0_{\tch}(g_u)=Z_G(g)/Z^0_G(g);$$
we denote by $Z'_G(g)$ the inverse image of $I_{\ch}(g_u)$ under the obvious surjective homomorphism 
$Z_G(g)@>>>Z_G(g)/Z^0_G(g)$ and by $Z'_{\tch}(g_u)$ the inverse image of 
$I_{\ch}(g_u)$ under the obvious surjective homomorphism 
$Z_{\tch}(g_u)@>>>Z_{\tch}(g_u)/Z^0_{\tch}(g_u)$, We have
$Z'_G(g)=Z'_{\tch}(g_u)$. This is a normal subgroup of $Z_G(g)=Z_{\tch}(g_u)$. 

For $g\in G^{sp}$ we set 
$$\bZ_G(g)=Z_G(g)/Z'_G(g),$$
a finite group. 

We show:

(a) {\it Assume that $g\in G^{sp}$. Then $g^q\in G^{sp}$ and $Z_G(g)=Z_G(g^q)$, $Z'_G(g)=Z'_G(g^q)$.}
\nl
Let $r\ge1$ be such that $g_s^{q^r}=g_s$. 
We have clearly $Z_G(g_s)\sub Z_G(g_s^q)\sub Z_G(g_s^{q^r})=Z_G(g_s^q)$
hence $Z_G(g_s)=Z_G(g_s^q)$. We have 
$Z_G(g_u)\sub Z_G(g_u^q)$; as in the proof of 2.1(a), $g_u,g_u^q$ are 
conjugate in $Z_G^0(g_s)$, hence they are also conjugate in $G$ and
$Z_G(g_u)=Z_G(g_u^q)$. Since $Z_G(g)=Z_G(g_s)\cap Z_G(g_u)$ and 
$Z_G(g^q)=Z_G(g_s^q)\cap Z_G(g_u^q)$,
we deduce that $Z_G(g)=Z_G(g^q)$, hence $Z^0_G(g)=Z^0_G(g^q)$ and 
$Z_G(g)/Z_G^0(g)=Z_G(g^q)/Z^0_G(g^q)$.
Let $\ch=Z_G^0(g_s)=Z_G^0(g_s^q)$. By 2.1(a), we have $I_{\ch}(g_u)=I_{\ch}(g_u^q)$. It
follows that $Z'_G(g)=Z'_G(g^q)$, as required.

We show:

(b) {\it Let $g\in G^{sp}$. Then $g_s^qg_u\in G^{sp}$ and $Z_G(g)=Z_G(g_s^qg_u)$,
$Z'_G(g)=Z'_G(g_s^qg_u)$.}
\nl
Let $\tch=Z_G(g_s)$. As in the proof of (a) we have $Z_G(g_s^q)=\tch$.
We have $Z_G(g)=Z_{\tch}(g_u)$ and $Z_G(g_s^qg_u)=Z_{\tch}(g_u)$; moreover,
$Z'_G(g)=Z'_{\tch}(g_u)$ and $Z'_G(g_s^qg_u)=Z'_{\tch}(g_u)$; (b) follows.

\mpb

The following result is immediate.

(c) {\it If $\fs:G@>>>G$ is an automorphism and $g\in G^{sp}$ then
$g'\m\fs(g')$ defines an isomorphism $Z'_G(g)@>>>Z'_G(\fs(g))$.}

\subhead 2.4\endsubhead
Let $g\in G$ and let $\tch=Z_G(g_s),\ch=\tch^0$. We show:

(a)  $j(g^q)$ is $G$-conjugate to $j(g^q_sg_u)$.
\nl
It is enough to show that
$g^q$ is $G$-conjugate to $g^q_sg_u$. It is also enough to show that
$g_u^q$ is $\ch$-conjugate to $g_u$. This holds since $p=0$.

\mpb

Let $G^{sp}_*$ be the set of all $g\in G^{sp}$ which satisfy:

(i) $g$ is $G$-conjugate to $j(g^q)$; $g$ is $G$-conjugate to $j(g_s^qg_u)$ (if $p=0$);
\nl
(The two conditions in (i) are equivalent by (a).)

Clearly, $G^{sp}_*$ is a union of $G$-conjugacy classes.

\subhead 2.5\endsubhead
Let $g\in G^{sp}_*$. For any $e\in\ZZ$ we set 
$$H_e(g)=\{x\in G;xgx\i=j^e(g^{q^e})\},\qua{}'H_e(g)=\{x\in G;xgx\i=j^e(g_s^{q^e}g_u)\}.$$
Note that $H_1(g)\ne\emp,{}'H_1(g)\ne\emp$. We set $H(g)=\sqc_{e\in\ZZ}H_e(g)$,
${}'H(g)=\sqc_{e\in\ZZ}{}'H_e(g)$. 
We regard $H(g)$ as a group with multiplication defined by 
$x*x'=j^{e'}(x)x'$ for $x\in H_e(g)$, $x'\in H_{e'}(g)$; we regard ${}'H(g)$ as a group with multiplication 
defined by $x*x'=j^{e'}(x)x'$ for $x\in{}'H_e(g)$, $x'\in{}'H_{e'}(g)$. 
Note that the map $H(g)@>>>\ZZ$ (resp. ${}'H_e(g)$) which takes $H_e(g)$ (resp. ${}'H_e(g)$) to $e$ 
for any $e\in\ZZ$ is a (surjective) group homomorphism with kernel $H_0(g)=Z_G(g)$
(resp. ${}'H_0(g)=Z_G(g)$). Thus, 
$Z'_G(g)$ can be regarded as a subgroup of $H(g)$ and also as a subgroup of ${}'H(g)$. It is in fact a 
normal subgroup. (It is enough to show that for $x_0\in Z'_G$ and $x$ in $H_e(g)$ or on ${}'H_e(g)$ we have
$x\i j^e(x_0)x\in Z'_G$. We have $j^e(x_0)\in Z'_G(j^e(g))$ hence 
$x\i j^e(x_0)x\in Z'_G(x\i j^e(g)x)$; by 2.3(a),(b), the last group equals 
$Z'_G(x\i j^e(g^q)x)=Z'_G(g)$ if $x\in H_e(g)$ and to 
$Z'_G(x\i j^e(g_s^qg_u)x)=Z'_G(g)$ if $x\in{}'H_e(g)$.)

We set $\D(g)=H(g)/Z'_G(g)$, ${}'\D(g)={}'H(g)/Z'_G(g)$. These have group structures induced from those
of $H(h),{}'H)g)$. 
The homomorphisms $H(g)@>>>\ZZ$, \lb ${}'H(g)@>>>\ZZ$ (as above) induces (surjective) group homomorphisms
$\et:\D(g)@>>>\ZZ$, ${}'\et:{}'\D(g)@>>>\ZZ$ with finite kernel $Z_G(g)/Z'_G(g)=\bZ_G(g)$. Thus 
$\et,{}'\et$ are $\ZZ$-gradings for $\D(g),{}'\D(g)$.

\subhead 2.6\endsubhead
Let $g\in G^{sp}_*$. Let $\ch=Z^0_G(g_s)$ and let
$\z:\l\m\z_\l$ be a homomorphism of algebraic groups 
$\kk^*@>>>\ch$ (coming from the Jacobson-Morozov theorem applied to 
$u\in\ch$) such that $\z_\l u\z_\l\i=u^\l$ for any $\l\in\kk^*$).
We define a map $\Ps:{}'H(g)@>>>H(g)$ by
$x\m y$ where $x\in{}'H_e(g),y=j^e(\z_{q^e})x, e\in\ZZ$. The computation
$$\align&ygy\i=j^e(\z_{q^e})xgx\i j^e(\z_{q^e}\i)=j^e(\z_{q^e})j^e(g_s^{q^e}g_u)
j^e(\z_{q^e}\i)\\&=j^e(\z_{q^e}g_s^{q^e}g_u\z_{q^e}\i)=j^e(g_s^{q^e}g_u^{q^e})
=j^e(g^{q^e})\endalign$$
shows that $\Ps$ is well defined and it maps ${}'H_e(g)$ to $H_e(g)$. It is in 
fact a bijection with inverse given by $y\m x$ where $y\in H_e(g)$,
$x=j^e(\z_{q^{-e}})y$, $e\in\ZZ$. 
For $x\in{}'H_e(g),x'\in{}'H_{e'}(g)$ we have

$\Ps(x*x')=j^{e+e'}(\z_{q^{e+e'}})j^{e'}(x)x'$,

$\Ps(x)*\Ps(x')=j^{e+e'}(\z_{q^e})j^{e'}(x\z_{q^{e'}})x'$
\nl
hence 
$$\Ps(x*x')=\Ps(x)*\Ps(x')c\tag a$$
where 
$$c=x'{}\i j^{e'}(x\z_{q^{e'}})\i j^{e+e'}(\z_{q^{e'}})j^{e'}(x)x'
=x'{}\i j^{e'}(\z_{q^{e'}}\i x\i j^e(\z_{q^{e'}})x)x'.$$
We show that 
$$c\in Z_G(g)^0.\tag b$$
We first show that for any $\l\in\CC$
$$c_\l:=x'{}\i j^{e'}(\z_\l\i x\i j^e(\z_\l)x)x'$$
belongs to $Z_G(g)$, or equivalently that
$$x'c_\l x'{}\i\in Z_G(x'gx'{}\i)=Z_G(j^{e'}(g_s^{q^{e'}}g_u))=j^{e'}(Z_G(g)),$$
It is enough to show that $\z_\l\i x\i j^e(\z_\l)x\in Z_G(g)$. We have 
$$\align&\z_\l\i x\i j^e(\z_\l)xg=\z_\l\i x\i j^e(\z_\l)j^e(g_s^{q^e}g_u)x\\&=
\z_\l\i x\i j^e(\z_\l g_s^{q^e}g_u)x=\z_\l\i x\i j^e(g_s^{q^e}g_u^\l z_\l)x\\&=
\z_\l\i g_s x\i j^e(g_u^\l z_\l)x=
\z_\l\i g_s g_u^\l x\i j^e(z_\l)x=gz_\l\i x\i j^e(z_\l)x\endalign$$
so that we have indeed $\z_\l\i x\i j^e(\z_\l)x\in Z_G(g)$. Thus we have
$c_\l\in Z_G(g)$ for any $\l\in\CC$. Note that $c_0=1$. Since $\l\m c_\l$ is a 
morphism $\CC@>>>Z_G(g)$ whose image contains $1$, its image must be contained in
$Z_G^0(g)$. Since $c=c_{q^{e'}}$, we see that (b) holds.
From (a),(b) we see that for $x,x'$ in ${}'H(g)$ we have
$$\Ps(x*x')\in\Ps(x)*\Ps(x')*Z'_G(g)\tag c$$
(recall that $Z^0_G(g)\sub Z'_G(g)$.) 
From the definitions we have $\Ps(x)=x$ for any $x\in{}'H_0(g)$; in particular,
$\Ps(Z'_G(g))=Z'_G(g)$. Using this and (c) 
we see that for $x\in{}'H(g)$ we have $\Ps(x*Z'_G(g))\sub\Ps(x)*Z'_G(g)$ so that
$\Ps$ induces a map $\bar\Ps:{}'\D(g)@>\si>>\D(g)$; here we view ${}'\D(g)$ (resp.
$\D(g)$) as the set of orbits of the right translation action of $Z'_G(g)$ on 
${}'H_(g)$ (resp. $H(g)$). From (c) we see that $\bar\Ps$ is a group homomorphism
(in fact a group isomorphism). Note also that $\bar\Ps$ induces the identity map
$\bZ_G(g)={}'\D(g)@>>>\D(g)=\bZ_G(g)$ and it is compatible with the $\ZZ$-gradings.

We now try to see to what extent does $\Ps$ depend on the choice of $\z$.
By a result of Kostant, any other choice of $\z$ must be of the form 
$\z':\l\m h\z_\l h\i$ for some $h\in Z^0_{\ch}(g_u)=Z^0_G(g)$. If we define 
$\Ps$ in terms of $\z'$ we obtain a bijection $\Ps':{}'H(g)@>>>H(g)$ given on 
${}'H_e(g)$ by
$$\Ps':x\m j^e(\z'_{q^e})x=j^e(h)\Ps(x)h\i hx\i j^e(h\i)x=h*\Ps(x)*h\i*d$$
where $d=hx\i j^e(h\i)x$. We show that
$$hx\i j^e(h\i)x\in Z_G^0(g)\text { for any }h\in Z^0_G(g).\tag d$$
We have 
$$\align&hx\i j^e(h\i)xg=hx\i j^e(h\i)j^e(g_s^{q^e}g_u)x=
hx\i j^e(h\i g_s^{q^e}g_u)x
\\&=hx\i j^e(g_s^{q^e}g_uh\i)x=hg x\i j^e(h\i)x=gh x\i j^e(h\i)x.\endalign$$
Thus, $hx\i j^e(h\i)x\in Z_G(g)$ for any $h\in Z^0_G(g)$.
We see that $h\m hx\i j^e(h\i)x$ is a morphism $Z^0_G(g)@>>>Z_G(g)$. Its image is
connected and its value at $h=1$ is $1$; hence its image is contained in $Z^0_G(g)$.
This proves (d).
From (d) we see that $\Ps$ and $\Ps'$ induce the same isomorphism
${}'\D(g)@>>>\D(g)$. Thus the isomorphism ${}'\D(g)@>>>\D(g)$ that we constructed is
independent of choices so it is canonical; we use it to identify these two groups.

\subhead 2.7\endsubhead
If $\D$ is a group with a given $\ZZ$-grading $\et:\D@>>>\ZZ$, we say that an automorphism $f:\D@>>>\D$
is inner if it is of the form $x\m x_0xx_0\i$ for some $x_0\in\et\i(0)$.
We consider an index set $I$ and a family $(\D_\a,\car_{\a\a'})$ where for any $\a\in I$, $\D_\a$ is 
a group with a given $\ZZ$-grading $\et$ and for any $\a,\a'$ in $I$, $\car_{\a\a'}$ is a nonempty set of
group isomorphisms $\D_{\a'}@>>>\D_\a$ compatible with the $\ZZ$-gradings.
We say that this family is coherent if the following conditions are satisfied.

(i) For any $\a\in I$, $\car_{\a\a}$ consists of the inner automorphisms of $\D_\a$.

(ii) For any $\a,\a'$ in $I$ and any $f\in\car_{\a,\a'}$ we have
$\car_{\a,\a'}=\car_{\a,\a}f=f\car_{\a',\a'}$.

(iii) For any $\a,\a',\a''$ in $I$, composition defines an inclusion 
$\car_{\a\a'}\car_{\a'\a''}\sub\car_{a,\a''}$.

In this case, for $i=0,1$ and for any $\a,\a'$ in $I$, the bijection
$M^i_{\et}(\D_{\a'})@>\si>>M^i_{\et}(\D_\a)$ (notation of 1.3)
induced by     $f:\D_{\a'}@>>>\D_\a$ (in $\car_{\a,\a'}$)
is independent of the choice of $f$. (We use the fact that an inner automorphism of $\D_\a,\et_\a$
induces the identity map on $M^i_{\et}(\D_\a)$.) It follows that the sets
$M^i_{\et}(\D_\a)$ with $\a\in I$ can be canonically identified with a single set associated to the
family.

\subhead 2.8\endsubhead
Let $\un G^{sp}$, $\un G^{sp}_*$ be the set of $G$-conjugacy classes contained in $G^{sp}$, $G^{sp}_*$. 
Let $C\in\un G^{sp}_*$. 
We consider the family $(\D(g),\car_{g,g'})$ indexed by $C$, where for $g\in C$, $\D(g)$ is the group with
a $\ZZ$-grading defined in 2.5 and for $g,g'$ in $C$, $\car_{g,g'}$ consists of the isomorphisms
$\t_k:\D(g')@>>>\D(g)$ (for various $k\in G$ such that $kg'k\i=g$) where
$\t_k$ is induced by the isomorphism $H(g')@>>>H(g)$ whose restriction $H_e(g')@>>>H_e(g)$ is 
$x\m j^e(k)xk\i$. From the definitions we see that this family is coherent. Hence for $i=0,1$, the sets
$M^i_{\et}(\D(g))$ with $g\in C$ can be canonically identified (as in 2.7)
with a single set $M^i(C)$ associated to $C$. We set 
$$AC(G)=\sqc_{C\in\un G^{sp}_*}M^0(C),\tag a$$
$$\cc(G)=\sqc_{C\in\un G^{sp}_*}M^1(C).\tag b$$
(Note that $AC(G),\cc(G)$ depend on $q$ and $j$.)
 The set $\cc(G)$ appeared in \cite{\ICM}.

\subhead 2.9\endsubhead
We define $\ph:G^{sp}@>>>G^{sp}$ by $g\m j(g^q)$. This is well defined by 2.3(a).
Clearly, $\ph$ induces a map $\un G^{sp}@>>>\un G^{sp}$ denoted again by $\ph$. We show:

(a) {\it Any orbit of $\ph:\un G^{sp}@>>>\un G^{sp}$ is finite. 
In particular, $\ph:\un G^{sp}@>>>\un G^{sp}$ is bijective.}
\nl
Let $g\in G^{sp}$. We choose $k\ge1$ such that $j^k=1$ and 
$g_s^{q^k}=g_s$. Then we have $\ph^k(g)=j^k(g^{q^k})=g_sg_u^{q^k}$ and 
this is conjugate to $g$. This proves (a). From the definitions we have

(b) {\it $\un G^{sp}_*$ 
is exactly the fixed point set of $\ph:\un G^{sp}@>>>\un G^{sp}$.}
\nl
We show:

(c) {\it $\ph:G^{sp}@>>>G^{sp}$ is a bijection.}
\nl
Let $g\in G^{sp},g'\in G^{sp}$ be such that $j(g^q)=j(g'{}^q)$. Then
$g^q=g'{}^q$ hence $g_s^q=g'_s{}^q$ and $g_u^q=g'_u{}^q$. The last 
equality implies $g_u=g'_u$ since $p=0$. From
$g_s^q=g'_s{}^q$ we deduce using 
$g_s^{q^r}=g_s$, $g'_s{}^{q^{r'}}=g'_s$ (where we can assume $r=r'\ge1$)
that $g_s=g'_s$ hence $g=g'$. Thus the map in (c) is injective.
Since $\ph:\un G^{sp}@>>>\un G^{sp}$ is surjective it follows that the
map in (c) is surjective. This proves (c).

\subhead 2.10\endsubhead
Let $C\in\un G^{sp}$. Let $\tC=\sqc_{g\in C}\bZ_G(g)$. We can view $\tC$ as an algebraic variety
which is a finite covering of $C$ under the obvious map $\io:\tC@>>>C$. Note
that $G$ acts on $\tC$ by $h:(g,zZ'_G(g))\m(hgh\i,hzh\i Z'_G(hgh\i))$;  
here $g\in G,z\in Z_G(g)$. This is well defined since $Z'_G(hgh\i)=hZ'_Gh\i$.
Now the orbits of $G$ on $\tC$ are exactly the connected components 
of $\tC$. Let $M(C)$ be the set of all pairs $(\co,\ce)$ where $\co$ is a
$G$-orbit in $\tC$ and $\ce$ is an irreducible $G$-equivariant local system
on $\co$. For any $g\in C$ there is a canonical bijection 
$M(C)@>\si>>M(\bZ_G(g))$ (notation of 1.1) which associates to $(\co,\ce)$ the equivalence 
class of $(x,\s)$ where $x$ belongs to $\co\cap\io\i(g)$ (a 
$\bZ_G(g)$-conjugacy class in $\bZ_G(g)$) and $\s$ is the stalk of $\ce$ at $x$
viewed as an irreducible representation of $Z_{\bZ_G(g)}(x)$. Thus the sets
$M(\bZ_G(g))$ with $g\in C$ are all in canonical bijection with $M(C)$. Let
$$CS(G)=\sqc_{C\in\un G^{sp}}M(C).$$

\subhead 2.11\endsubhead
For any $g\in G^{sp}$ we define a homomorphism $\bZ_G(g)@>>>\bZ_G(\ph(g))$
by \lb $zZ'_G(g))\m j(z)Z'_G(j(g^q))$ with $z\in Z_G(g)$. We show:

(a) {\it this map is a bijection.}
\nl
Since $\bZ_G(g),\bZ_G(\ph(g))$ are finite, it is enough to show that our map 
is injective. 
Let $z,z'$ in $Z_G(g)$ be such that
$j(z)=j(z')\mod Z'_G(j(g^q))$. Then $z=z'\mod Z'_G(g^q)=Z'_G(g)$, see 
2.3(a). This proves injectivity hence (a).

\mpb

We define $\ti\ph:\sqc_{C\in\un G^{sp}}\tC@>>>\sqc_{C\in\un G^{sp}}\tC$
by $(g,zZ'_G(g))\m(j(g^q),j(z)Z'_G(j(g^q)))$; here $g\in C,z\in Z_G(g)$.
Using (a) and 2.9(a) we see that $\ti\ph$ is a bijection.

For any $C\in\un G^{sp}$, $\ti\ph$ restricts to a (bijective) morphism
$\tC@>>>\tC'$ where $C'=\ph(C)$.
Note that the $G$-action on $\tC$ is compatible under $\ti\ph$ with the
$G$-action on $\tC'$ composed with the automorphism $j$ of $G$; in 
particular, $\ti\ph$ carries any $G$-orbit on $\tC$ to a $G$-orbit on 
$\tC'$. We see that $(\co,\ce)\m(\ti\ph(\co),\ti\ph_!\ce)$ defines a 
bijection $M(C)@>\si>>M(C')$. Taking disjoint unions we see that $\ti\ph$
induces a bijection
$$CS(G)@>\si>>CS(G).\tag b$$
The fixed points of this bijection are the pairs $(\co,\ce)\in M(C)$
with $\ph(C)=C$ such that $\co=\ti\ph(\co)$, $\ce\cong\ti\ph_!\ce$.
The condition that $\ph(C)=C$ is the same as $C\in G_{j,q}$, see 2.9(b).
Assuming that this is the case, let $g\in C$. Then $\D(g),\et$ and
$M^1_\et(\D(g))$ are well defined and $\bZ_G(g)$, $M(\bZ_G(g))$ are well
defined. Note that $\bZ_G(g)$ is the subgroup $\et\i(0)$ of $\D(g)$.
Hence $M^1_\et(\D(g))$ can be viewed as the subset of $M(\bZ_G(g))=M(C)$
defined by the two conditions in 1.3(c). In terms of $M(C)$ those two
conditions are exactly the conditions that $\co=\ti\ph(\co)$, 
$\ce\cong\ti\ph_!\ce$ (as above). In this way we see that 

(c) {\it $AC(G)$ (see 2.8) can be naturally identified with the
fixed point set of the map (b) of $CS(G)$ into itself.}

\subhead 2.12\endsubhead
In the reminder of this section we assume that $q$ is a power of a
prime number and that $\cg$ is a connected reductive group over $\bar\FF_q$
of type dual to that of $G$ with a fixed $\FF_q$-rational structure which is
twisted according to $j$. Let $\cf$ be the vector space of class
functions $\cg(\FF_q)@>>>\CC$. 
Now $\cf$ has a basis $\b$ given by the characters of the irreducible
representations of $\cg(\FF_q)$. We define a function $\d:\b@>>>\{1,-1\}$ as
follows: if $\c$ is the character of an irreducible representation $\e$ of 
$\cg(\FF_q)$ and $D(\e)$ is its dual (as in \cite{\DL}) then $\d(\c)$ is
characterized by the property that $\d(\c)D(\e)$ is an irreducible representation.

Another basis $\b'$ of $\cf$
(whose elements are defined only up to multiplication by roots of $1$) is 
provided by the characteristic functions of the character sheaves
on $\cg$ which are isomorphic to their inverse image under the Frobenius map
$\cg@>>>\cg$. It is known that $\b$ is indexed by $\cc(G)$, the character
sheaves on $\cg$ are indexed by $CS(G)$ and the Frobenius-invariant character 
sheaves on $\cg$ are indexed by the fixed point set of $\ph$ on $CS(G)$. Thus, 
$\b'$ is indexed by $AC(G)$. 

The decompositions 2.8(a),(b) give rise to decompositions
$\b=\sqc_C\b_C,\b'=\sqc_C\b'_C$. 
We now fix $C\in\un G^{sp}_*,g\in C$ and we list the elements of $\b_C$ as 
$\{\bb_{x,\s};(x,\s)\in M^1_\et(\D)\}$ and the elements of $\b'_C$ as 
$\{\bb'_{x,\s};(x,\s)\in M^0_\et(\D)\}$ (where $\D=\D(g)$).

For any $(x,\s)\in M^0_\et(\D)$ we set
$$\align&\bb''_{x,\s}\\&=\sum_{(y,\t)\in M^1_\et(\D)}
\sum_{w\in\D/\la\x\ra;xwyw)\i=wyw\i x}\et(\x)
|Z_\D(x)/\la\x\ra|\i|Z_\D(y)/\la\x\ra|\i\\&
\T\tr(w\i x\i w,\ti\t)\tr(wyw\i,\ti\s)\d(\bb_{y,\t})\bb_{y,\t}\tag a\endalign$$
where  $\x\in\D^*$ (notation of 1.3),
$\ti\s$ is a representation of $Z_\D(x)$ extending $\s$ on which $\x$ acts as $1$,
$\ti\t$ is the representation of $Z_\D(y)$ extending $\t$ on which
$y$ acts as $1$. Note that $\bb''_{x,\s}$ is independent (up to 
multiplication by a root of $1$) of the choice of $\x,\ti\s$, see 1.3.

We expect that for any $(x,\s)\in M^0_\et(\D)$ we have
$$\bb'_{x,\s}=\bb''_{x,\s}\tag b$$
up to multiplication by a root of $1$.

Note that (b) is known in many cases but not in complete generality.
It implies that:

(c) {\it for any $C\in\un G^{sp}_*$, the subspace of $\cf$ spanned by $\b_C$ 
coincides with the subspace of $\cf$ spanned by $\b'_C$.}
\nl
The functions $b''_{x,\s}$ are called {\it almost characters}. Note that the 
almost characters (for various $C$) form a basis of $\cf$ indexed by $AC(G)$.


\widestnumber\key {AB}
\Refs
\ref\key{\DL}\by P.Deligne and G.Lusztig\paper Duality for representations
of reductive groups over a finite field\jour J.Alg\vol74\yr1982\pages284-291\endref
\ref\key{\EEIGHT}\by G.Lusztig\paper Unipotent representations of a finite 
Chevalley group of type $E_8$\jour Quart.J.Math.\vol 30\yr1979\pages 315-338\endref
\ref\key{\ICM}\by G.Lusztig\paper Characters of reductive groups over 
finite fields\inbook Proc. Int. Congr. Math. Warsaw 1983\publ
 North Holland\yr1984\pages877-880\endref
\ref\key{\ORA}\by G.Lusztig\book Characters of reductive groups over a finite 
field \bookinfo Ann. Math.Studies 107\publ Princeton U.Press \yr1984\endref 
\ref\key{\FAM}\by G.Lusztig\paper Families and Springer's correspondence\jour 
Pacific J.Math.\vol 267\yr2014\pages 431-450\endref
\endRefs
\enddocument